\input amstex
\magnification =1200
\documentstyle {amsppt}

\define\mas{Monge-Amp\`ere solution}

\define\qh{quasi-homogeneous}

\define\gt{Grauert tube}
\define\Rm{Riemannian manifold}
\document
\baselineskip 1.4\baselineskip
\topmatter
\title  On rigidity of Grauert tubes \endtitle
\author Su-Jen Kan 
\endauthor 
\address 
Institute of Mathematics, Academia Sinica, Taipei 115, Taiwan
\endaddress
\email 
{kan\@math.sinica.edu.tw}
\endemail
\date Oct. 8, 2003
\enddate
\NoRunningHeads
\thanks 2000 {\it Mathematics Subject Classification. } 32C09, 32Q28, 32Q45, 53C24, 58D19. 
\endthanks
\abstract

Given a real-analytic Riemannian manifold $M$ there exists a canonical  complex structure on part of its
tangent bundle which turns leaves of the
Riemannian foliation on
$TM$ into holomorphic curves. A {\it Grauert tube} over $M$ of radius $r$,
denoted as $T^rM$, is the collection of tangent vectors of $M$ of length less than $r$ equipped
with this canonical complex structure. We say the Grauert tube $T^rM$ is {\it rigid} if $Aut(T^rM)$ is
coming from $Isom (M)$.

In this article, we prove the  rigidity for Grauert tubes over quasi-homogeneous
Riemannian manifolds. A Riemannian manifold
$(M,g)$ is quasi-homogeneous if  the quotient space
$M/Isom_0 (M)$ is compact. This category has included compact \Rm s, homogeneous Riemannian
manifolds, co-compact \Rm s whose isometry groups have dimensions $\ge 1$, and products of the above
spaces.

\endabstract
\endtopmatter

\subheading {1. Introduction}

The {\it adapted complex
structure} is the unique complex structure on (part of) the tangent bundle of a real-analytic
Riemannian manifold $(M,g)$ which turns  leaves of the Riemannian foliation on $TM$ into 
holomorphic curves. It was proved by Guillemin-Stenzel and independently by Lempert-Sz\H oke that
such an adapted complex structure exists  in a neighborhood of $M$ in $TM$ by solving a complex
homogeneous Monge-Amp\`ere equation. We make a
remark here that the \mas s $\rho_{\sssize GS}$ in [G-L],  $\rho_{\sssize LS}$ in [L-S] and
$\rho_{\sssize K}$ in [K1] all have different scalings. The relation among these three are: 
$\rho_{\sssize GS}=2 \rho_{\sssize K}=4 \rho_{\sssize LS}.$ The $\rho_{\sssize LS}$ is exactly the
length square function. Therefore, we will adapt the scaling that used in [L-S]; the $\rho$ we
consider in this article is the \mas\ with the initial condition  $\rho_{i\bar j}|_
{\sssize M}=\frac 12 g_{ij}$ and $\rho (x,v)=|v|_g^2$ for any tangent vector $v\in T_xM.$
 The set of tangent vectors of length less than
$r$ equipped with the adapted complex structure is called a {\it Grauert tube } of radius $r$ centered
at $M$. It is denoted as 
$T^rM$. In terms of the \mas, $T^rM=\{\rho^{-1}[0,r^2)\}.$ 

Associated to each $M$ there  exists a maximal possible radius $r_{max}(M)\ge0$ such that $T^rM$ exists
for any $r\le r_{max}(M).$ Since the adapted complex structure is
locally defined, the maximal possible radius $r_{max}(M)>0$ for compact $M$.  Nevertheless,
$r_{max}(M)$ might very well be zero for non-compact $M$. 

Since the
adapted complex structure is constructed canonically associated to the Riemannian metric
$g$ of the  manifold $M$, the differentials of the isometries of $(M,g)$ are
natural automorphisms of $T^rM$. On the other hand, it is interesting to see whether every automorphism
of $T^rM$ comes from this way or not. When the answer is
affirmative, we say the \gt\ is {\it rigid}. 

With respect to the adapted complex structure, the length square function $\rho$ is strictly
plurisubharmonic. When the center $M$ is compact the
\gt \ $T^rM$ is Stein since it is   exhausted by the strictly plurisubharmonic function $\rho$.
Furthermore, if the radius $r$ is less than the maximal possible radius, $T^rM$ is then a bounded
domain with smooth
strictly pseudoconvex boundary in a Stein manifold. In this case, most good properties of  bounded
strictly pseudoconvex domains are inherited  and the
automorphism group of $T^rM$ is a compact Lie group. Base on these, Burns and Hind ([B-H]) are
able to prove the rigidity for \gt s over compact real-analytic Riemannian manifolds. 

For non-compact $M$, the length square function $\rho$ is no longer an exhaustion. This
has made the situation complicated to deal with. In fact, nothing particular is known, not even to
the  existence of a
\gt\ over non-compact $M$. When the \gt s exist, most of the good properties in the compact
case are lacking here since the length square function $\rho$  no longer  exhausts the tube.

In [K2], using the homogeneity of the central manifolds, the author was able to prove 
the rigidity of Grauert tubes over homogeneous Riemannian manifolds. The author also 
 proved a  weak rigidity for general Grauert
tubes which says that if  $r_{max}(M)>0$ and $T^rM$ is not covered by
the ball,  then
$Aut_0(T^rM)\approxeq Isom_0(M)$ for any
$r<r_{max}(M)$.

It does not seem clear to us what kind of geometric properties possessed by the real-analytic Riemannian
manifold
$(M,g)$ will guarantee  the existence the Grauert
tube structure, i.e., will enforce $r_{max}(M)>0.$ To the author's knowledge, there are only four
categories that we know:  $M$ is compact, is homogeneous, is co-compact or is the product of the above
cases.  The rigidity of Grauert tubes over the first two cases have been clarified in [B-H] and [K2],
respectively. The motivation of this article is to examine the rigidity of \gt s over
the third  and the fourth cases.  

In
an attempt to tackle this problem, we formulate the situation as below.
 The objects we are interested in are real-analytic Riemannian manifolds $M$ possesses the property that
every point in $M$ could be sent to a compact subset of $M$ through some isometry.
  First of all,  the above four
cases(compact, homogeneous, co-compact and products of them) are all included in this family.
Secondly, this is a kind of  generalized homogeneous spaces; lots of  the technique in the
homogeneous case could be transplanted here without much difficulty. However, to prove the rigidity, we
need a reasonable large isometry group. So we defined a \Rm\  $M$ is  quasi-homogeneous if
$Isom_0(M)\cdot K=M$ for some compact subset $K$ of $M$. This would immediately implies that $ dim
Isom_0(M)\ge 1$ for non-compact $M$.
The main result of this
article is Theorem 4.5 which proved the rigidity of Grauert tubes  over \qh\ \Rm s. 
 \proclaim{Theorem 4.5}
Let $M$ be  a \qh\ manifold. Let
$T^rM$ be the
\gt\ of radius
$r<r_{max}(M)$ which is not
 covered by the ball.  Then
$Isom(M)\approxeq Aut(T^rM)$ and $T^rM$ has a unique center.
\endproclaim

The organization of this article is as the following. In $\S2$ we define quasi-homogeneity and derive
some basic properties of quasi-homogeneous  manifolds. We prove that a Grauert tube over
quasi-homogeneous  manifold is complete hyperbolic in $\S3$. The $\S4$ is devoted to the proof of
the rigidity of $T^rM$ for quasi-homogeneous $M$.

The author would like to thank Professor Kang-Tae Kim for bringing her attention to the case
$M/Isom(M)$ is compact.

\subheading{\S 2 Quasi-homogeneous manifolds}

Recall  that a \Rm \ $M$ is called homogeneous if for any two points $x,y\in M$, there exists a $g\in
Isom(M)$ such that $g\cdot x=y$. This transitivity property has implied that the isometry group of $M$
is very large, it is at least as large as $M$.

Motivated by the homogeneous case, the real-analytic Riemannian manifolds  we consider through out this
article are Riemannian  manifolds possess certain kind of transitivity and have
 reasonably large
isometry groups.

\definition{Definition}
A real-analytic Riemannian manifold $M$ is  {\it\qh } if  the quotient space $M/Isom_0(M)$ is compact.
\enddefinition
The condition is equivalent to  that there exists a compact subset $K$ of
$M$ such that for any
$m\in M$ there associated   a $g\in Isom_0 (M)$ with $g\cdot m\in K$ or equivalently, $Isom_0(M)\cdot
K=M.$ We may further assume that both $M$ and $K$ are connected and  that for any $x\ne y, x,y \in
\overset {\sssize\circ }\to K$, there is no
$g\in Isom_0 (M)$ such that $g\cdot x=y.$

The reason we ask for $Isom_0 (M)$ instead of $Isom (M)$ is that the criteria would force the isometry
group of $M$ to be fairly large for  non-compact manifold $M$, say $\dim
Isom (M)\ge 1$. In fact, if a non-compact connected manifold $M$ has a connected compact quotient
 $M/Isom(M)$ then either $\dim
Isom (M)\ge 1$ which would automatically imply that $M/Isom_0(M)$ is compact   or $M$ is a
non-compact co-compact manifold with a discrete isometry group. 

 This category actually has
contained most of cases that we are sure about the existence of Grauert tubes, examples of such
$M$ are: compact manifolds,  homogeneous manifolds, co-compact manifolds whose isometry group have
dimensions $\ge 1$, and products of the above manifolds. The only thing out of control are non-compact
co-compact manifolds with discrete isometry groups. We examine some of the fundamental properties
of \qh\ spaces. 
\proclaim{Remark}
The properties proved in this and the next sections hold for any \Rm\  $M$ which possesses the property
that  there exists a compact subset $K\subset M$ with the property that  $Isom(M)\cdot K=M$. That is,
they also work for non-compact co-compact manifolds with discrete isometry groups. However, for the
rigidity proof in the last section of this article, we do need some dimension criteria on the isometry
group.
\endproclaim

\proclaim{Lemma 2.1}
If $(M,g)$ is a  \qh\ manifold, then  $(M,g)$ is complete.
\endproclaim 
\demo{Proof}To prove $(M,g)$ is complete is equivalent to show that every geodesic could be extended
to a geodesic defined on the whole real line. 

Let $K$ be a compact subset of $M$ such that every point of $M$ could be sent to $K$ by an 
isometry $h\in Isom_0 (M)$.
Let $s$ be a positive number such that the exponential map is defined on the ball $B(p,s)$ for all
$p\in K$. By hypothesis, any point $m$ in $M$ could be sent to a point in $K$ by an isometry.
Therefore, for any
$m\in M$, the exponential map $\exp_m$ is defined on $B(m,s)$. 

Let $\eta(t)$ be a geodesic with $\eta (0)=m$ and $I=(-a,b)$ be its maximal interval of 
definition. Then both $a$ and $b$ are greater than $s$.
Suppose
$b<\infty$, we set the point $p=\eta (b-\frac s2)$. Let $\gamma$ be the geodesic with $\gamma (0)=p$
and $\gamma'(0)=\eta'(b-\frac s2)$. Then this $\gamma$ extends $\eta$ and is at least defined
on the interval $(-s,s)$. That is, $\eta$  is defined on $(b-\frac {3s}2, b+\frac s2)$, a
contradiction. Similarly, we may show that $a=\infty$.
\qed
\enddemo
\proclaim{Lemma 2.2}
If $M$ is a  \qh\ manifold. Then  $r_{max}(M)>0$.
\endproclaim 
\demo{Proof}
Let $K\subset M$ be the compact set described as before. As the adapted complex structure is a
local object, there exists a $r>0$ such that the adapted complex structure is defined on
$T_p^rM:=\{(p,v):v\in T_pM, |v|<r\}$ for all $p\in K$. Since any point $m\in M$ could be sent to a point
in
$K$ by an isometry(hence an automorphism), the adapted complex structure is well-defined on $T_m^rM$
for all
$m\in M$. The lemma is claimed.\qed

\enddemo
\proclaim{Lemma 2.3}
Let $(M,g)$ be a non-compact \qh\ manifold and let
$K\subset M$ be a  compact subset described as above. Then for every 
 $p\in K$,  there exist
infinitely many
$m\in M$ and
$h_m\in Isom_0(M)$ such that $h_m(m)=p$.

\endproclaim 
\demo{Proof} It is clear that $Isom_0(M)$ is non-compact.
Since $K$ is compact in a complete Riemannian manifold, it is bounded. We may pick a large
geodesic ball
$B:=B(p,s)$ centered at
$p$ of radius
$s$ containing
$K.$  Given $q_1\in M-K$, there exists $h_{1}\in Isom_0(M)$ such that $h_{1}(q_1)\in K$. We
may assume that $h_{1}(q_1)=p$. If not, there exists a large set $L$ outside of $K$ such that
$h_{1}(L)\supset K$, then we are able to find some $\hat q\in L$ such that $h_{1}(\hat q)=p$.

Denote $h_{1}^{-1}(K):=\hat K_1\subset \hat B_1:= h_{1}^{-1}( B)$, a geodesic ball of radius $s$
centered at $q_1$. Take $q_2\in M-\hat B_1$ such that $d(q_1,q_2)>2s$ and such that there exists a
$h_2\in Isom_0(M)$ sending $q_2$ to $p$ where $d$ denotes the distance function induced by the
Riemannian metric $g$. Denoting
$h_{2}^{-1}(K):=\hat K_2\subset
\hat B_2:= h_{2}^{-1}( B),$ then $\hat B_2\cap \hat B_1=\emptyset$.

As $M$ is complete, we are able to continue the process and find infinitely many different $q_j\in M$
and $h_j\in Isom_0(M)$ such that $h_j(q_j)=p, \forall j$.\qed
\enddemo

\

\subheading{\S 3 Complete hyperbolicity of $T^rM$}

The Kobayashi pseudo-metric is defined in any complex manifold. When it is a metric, we call such a
complex manifold a hyperbolic manifold. It was shown by Sibony in [S] that a complex manifold admitting
a bounded strictly plurisubharmonic function is  hyperbolic. In a \gt\ of finite radius, the
length square function $\rho$ is  strictly plurisubharmonic and bounded in the whole tube. Therefore,
every \gt\ of finite radius is hyperbolic. On the other hand, a \gt\ of infinite radius, i.e., the
adapted complex structure is defined on the whole tangent bundle of $M$, is never hyperbolic.
\proclaim{Proposition 3.1}
Let $TM$ be a \gt\ of infinite radius. Then $TM$ can not be hyperbolic. 
\endproclaim
\demo{Proof}
Following the definition of the adapted complex structure, for a given arc-length parametrized geodesic
$\gamma$ in
$M$ there exists a holomorphic mapping
$$ f: \Bbb C \to TM,\;\;\;
f(\sigma+i\tau )= \tau\gamma'(\sigma).$$ 
 Since
the Kobayashi metric in $\Bbb C$ is trivial, by Picard's
theorem, the mapping $f$ has to be constant if  $TM$ is hyperbolic. A contradiction.\qed

\enddemo
 In general, it is hard to see whether a hyperbolic  manifold  is
complete or not. If $M$ is compact and the radius $r<r_{max}(M)$ then the \gt\ $T^rM$ is complete
hyperbolic since it is a bounded domain in the
Stein manifold $T^{r_{max}}M$   with smooth strictly pseudoconvex boundary. The same holds for
co-compact $M$ since  $T^r (M/\Gamma)=T^rM/\Gamma$ for
$\Gamma<Isom(M)$ and the \gt\ $T^rM$ is the  covering of a complete hyperbolic manifold. 
 In [K2], we also
proved that $T^rM$ is complete hyperbolic when  $M$ is homogeneous and $r<r_{max}(M)$. 

In
[F-S] Fornaess and Sibony has proved a sufficient condition for a hyperbolic manifold to be complete
hyperbolic.  They show that   a hyperbolic manifold
$\Omega$ is complete hyperbolic if the quotient
$\Omega/Aut(\Omega)$ is compact. 
Inspired from Fornaess-Sibony's work, we  prove in this  section that the
Grauert tube
$T^rM$ is complete hyperbolic if $M$ is \qh\ and $r<r_{max}(M)$.

Let $d_{\sssize K}$ be the Kobayashi metric of the hyperbolic manifold $T^rM$ and
$\hat d_{\sssize K}$ be the restriction of $d_{\sssize K}$ to $M$. That is,  the metric 
$\hat d_{\sssize K}$ is defined as $$\hat d_{\sssize K}(p,q):=d_{\sssize K}(p,q), \forall p,q\in M.$$
From the construction of \gt s, the $Isom(M)$ is naturally included in the $Aut(T^rM)$. Therefore any
$h\in Isom(M)$ is an isometry for the metric $\hat d_{\sssize K}$.

Recall that a metric space  is complete if  every Cauchy sequence
converges. Since the Kobayashi metric $d_{\sssize K}$ is an inner metric (i.e., one that  comes from
the arc length), the Hopf-Rinow-Myers theorem shows that the completeness  is equivalent to that every
finite ball is relatively compact.  We will prove the completeness of $d_{\sssize K}$ by first showing
that
$(M,\hat d_{\sssize K})$ is complete, i.e., every Cauchy sequence converges.
\proclaim{Lemma 3.2}
$\hat d_{\sssize K}$ is a complete metric in $M$.
\endproclaim
\demo{Proof}
  We would like to show that every Cauchy sequence in $(M,\hat d_{\sssize K}) $  converges.
Since $M$ is \qh, there is a compact set $K\subset M$  such that for any 
$p\in M$ there exists a $h\in Isom_0(M)$ with $h(p)\in K$. Notice that by the definition of a
Riemannian manifold, the natural topology in
$M$ and the topology  induced by the metric $g$ are the same.

Since $K$ is compact in the complete Riemannian manifold $(M,g)$, $K$ is bounded and closed. We may 
find a ball $B_g(x_0,R)$, centered at $x_0\in M$ of radius $R$ with respect to the
$g$-metric, in $M$  such that
$K\subset B_g(x_0,R)$. Let
$F:=B_{g}(x_0, 2R)$. Since $d_{\sssize K}$ is continuous, the metric $\hat d_{\sssize K}$ is continuous
as well.  Therefore, there exists an $\epsilon >0$ such that $\hat d_{\sssize K}(p,q)\ge\epsilon$  for
any
$p\in K, q\notin \bar F$.

Let $\{p_j\}\subset M$ be a Cauchy
sequence in
$(M, \hat d_{\sssize K})$. That is, given $\epsilon >0$ there exists a
large
$m$ such that
$\hat d_{\sssize K}(p_k,p_l)<\epsilon,
\forall k, l\ge m.$ Let $h\in Isom_0(M)$ such that $h(p_m)\in K$. Then
$$\hat d_{\sssize K}(h(p_j),h(p_m))=d_{\sssize K}(h(p_j),h(p_m))=d_{\sssize
K}(p_j,p_m)=\hat d_{\sssize K}(p_j,p_m)<\epsilon,\forall j>m.$$
The Cauchy sequence $\{h(p_j)|j=m,\cdots, \infty\}$ lies in the compact set $\bar  F$. It converges
to a point $\hat p\in  \bar F$ and hence the Cauchy sequence $\{p_j\}$ converges to the point
$h^{-1}(\hat q)\in M$.
\qed
\enddemo

\proclaim{Theorem 3.3}

Let $M$ be a \qh\ manifold.  Then for any $r<
r_{max}(M)$,
$T^rM$ is complete hyperbolic.
\endproclaim 
\demo{Proof} 
Fix $z=(p_1,v_1)\in T^rM$. We compute the Kobayashi distance from $z$ to the point $w=(p_2,v_2)$.
Let $h_1\in Isom_0(M)$ such that $h_1(p_2)\in K$, then there exists a constant $L_1(v_2)$ such that 
$$d_{\sssize
K}((p_2,0),(p_2,v_2))=d_{\sssize K}((h_1(p_2),0),(h_1(p_2),h_{1*}v_2)\le L_1(v_2).\tag 3.1$$  Now,
$$\aligned
 d_{\sssize K}(z,w)&\ge d_{\sssize K}((p_1,0), (p_2,v_2))-d_{\sssize
K}((p_1,0), (p_1,v_1))\\
&\ge d_{\sssize K}((p_1,0), (p_2,0))-d_{\sssize
K}((p_2,0), (p_2,v_2))-d_{\sssize
K}((p_1,0), (p_1,v_1))\\
&=\hat d_{\sssize K}(p_1,p_2)-d_{\sssize
K}((p_2,0), (p_2,v_2))-d_{\sssize
K}((p_1,0), (p_1,v_1)). \endaligned\tag 3.2$$
We claim the complete hyperbolicity by showing that for any radius $R$ the  bounded ball $B_{\sssize K
}(z,R)$ is relatively compact in $T^rM$.

Since $T^rM$ is a domain in the complex manifold $T^{r_{max}}M$ and every point $(x,v)\in T_xM, |v|=r$
is a smooth strictly pseudoconvex boundary point of  $T^rM$. By Lemma 2.2 of [K2], if $w=(p_2,v_2)\in
B_{\sssize K }(z,R)$ then there exists $\delta>0$ such that $|v_2|<r-\delta$. By equations (3.1) and
(3.2),
$$\aligned
\hat d_{\sssize K}(p_1,p_2)&\le d_{\sssize K}(z,w) +d_{\sssize K}((p_2,0),
(p_2,v_2))+d_{\sssize K}((p_1,0), (p_1,v_1))\\
&<R+L_1(v_2)+ L_1(v_1)<L
\endaligned$$ for some positive $L$. Thus, by Lemma 3.2, $p_2$
lies in some bounded set in $M$. Thus $B_{\sssize K }(z,R)$ is relatively compact in $T^rM$ and the
Kobayashi metric is complete.
\qed
\enddemo

\

\subheading { 4.  The rigidity of \gt s} 

For compact $X$, Burns-Hind [B-H] have proved that the isometry group   $Isom (X)$ of $X$   is
isomorphic to the automorphism group $ Aut (T^rX)$ of the \gt\   for any radius $r\le r_{max}(X)$.

For the non-compact cases, the best rigidity results so far are in [K2]. It shows that the rigidity
holds for  any
\gt\ over a homogeneous \Rm\  of $r< r_{max}$; it also claims   that the identity component of the
automorphism group of the \gt\ is isomorphic to the identity component of the isometry group of the
center manifold for  \gt s over general real-analytic \Rm s 
 of $r< r_{max}$. The only exception for the above two results occurs  when the \gt\
 is covered by the ball. It was proved by the author in [K1] that  $T^rX$ is biholomorphic
to
$B^n\subset\Bbb C^n$ if and only if  $X$ is the real hyperbolic space $\Bbb H^n$ of curvature $-1$
and the radius
$r=\frac {\pi}{4}$. Apparently, the automorphism group of $B^n$ is much larger than the
isometry group of $\Bbb H^n$. We also remark here that the restriction $r<r_{max}$ is necessary as the
rigidity fails for \gt s over non-compact symmetric spaces of rank-one of maximal radius shown in
[B-H-H]. 

Let $\tilde M$ be the universal covering of $M$. We have proved in Lemma 4.2 of [K2] that if there
is a unique Grauert tube representation for the complex manifold $T^r\tilde M$ then there is a
unique Grauert tube representation for  $T^r M$. Denote
$I=Isom (M)$ and $G=Aut(T^rM)$; $I_0$ and $G_0$  the corresponding identity components. We have shown
in Theorem 6.4 of [K2] that $G_0=I_0$ provided that $T^rM$ is not covered by the ball.  We may assume
from now on  that $M$ is a simply-connected \qh\ manifold, $r<r_{max}(M)$ and that $T^rM$ is not covered
by the ball.   Since  $M$ is a connected \qh\ \Rm, there
exists a connected compact subset $K\subset M$ such that $I_0\cdot K=M$ and  for any  $x\ne
y, x,y\in \overset {\sssize\circ }\to K$, there is no
$g\in Isom_0 (M)$ such that $g\cdot x=y.$

For any given $f\in G$, the set $N=f(M)$ equipped with the push-forward metric coming from $M$
is another center of the Grauert tube $T^rM$. It follows from the Theorem 6.4 of [K2]
 that $Isom_0(N)=G_0=Isom_0(M)$. Therefore,
$I_0\cdot N=N.$

 Let's denote the projection as  $$ \pi: T^rM\to  M,
\;\;
\pi (p,v)=p.$$ It is clear that $\pi$ is
$I_0$-invariant since for any $(x,v)\in T^rM$ and $g\in I_0$, we have
$$\pi (g\cdot (x,v))=\pi (g\cdot x, g_* v)=g\cdot x=g\cdot \pi (x,v).\tag 4.1$$
The Riemannian metric $g$ of $M$ determines the Levi-Civit\`a connection on $TM$ which splits the
tangent space $T_z(TM), z\in TM,$ into vertical and horizontal spaces. The connection map is
$K:T_z(TM)\to T_z(T_{\pi (z)}M).$ A vector $\zeta\in T_z(TM)$ is horizontal  if $K\zeta=0$; is vertical
if
$\pi_*\zeta=0$.
Let's make some preliminary observation on the projection of $N\subset T^rM$.
\proclaim{Lemma 4.1} $\dim (\pi (N))=\dim M=n.$
\endproclaim
\demo{Proof}
Let $U_z\subset N$ be a small neighborhood of $z\in N$. We would like to show that its projection $\pi
(U_z)$ has dimension $n$. Suppose not, $\dim \pi (U_z)=\mu<n.$ Let $\dim I_0=l, \dim K=k.$  
Then the tangent space $T_z(U_z)$ is
spanned by $\mu$ horizontal vectors and $n-\mu$ vertical vectors $\{\eta_1,\cdots,\eta_{n-\mu}\}$. As
$N$ is $I_0$ invariant, and 
$I_0$ moves vertical vectors to vertical vectors. By the quasi-homogeneity of $M$, the dimension of
$N$ would be
$n-\mu + l+\mu$ if $\mu<k$; would be $n-\mu +n$ if $\mu\ge k$.  In either case,  the  dimension of $N$
would be greater $ n$ since
$l\ge 1$ and $n-\mu >0$. A contradiction. Therefore, $\mu=n$.\qed
\enddemo

We would like to show that  $N$ actually crosses through each fiber $T_p^rM:=\{(p,v): v\in T_pM,
|v|<r\}, p\in M$. That is, $\pi (N)=M$. 

\proclaim{Lemma 4.2} 
$N\cap T_p^rM\ne\emptyset$ for any $p\in M$.
\endproclaim
\demo{Proof} It is clear from the above lemma that $\pi (N)$ and $M$ have the same dimension and $\pi
(N)$ is connected since $\pi$ is a continuous map. It is therefore sufficient to claim that $\pi (N)$
is both open and closed. The connectedness of $M$ would immediately implies that $\pi (N)=M.$

For any $p\in \pi (N)$, there exists a $z\in N$ such that $\pi (z)=p$. Since $N$ is a submanifold,
every point is an interior point and there is a neighborhood $U_z$  of $z$ in $N$. By the proof of
Lemma 4.1, 
$\pi (U_z)$ has dimension $n$ and hence  is a neighborhood of $p$ in $\pi (N)$. Therefore, $\pi (N)$ is
open in
$M$.

For the closeness. We first observe that there exists a real number $0<l<r$ such that every $(x,v)\in
N$ has $|v|<l$. This comes from the fact that $I_0\cdot K=M$ and $N=f(M)=f(I_0\cdot K)=I_0\cdot f(K).$
 Let $\{x_j\}$ be a sequence in $\pi (N)$ with $\lim_{j\to\infty} x_j=x\in M$. Then for give
$\epsilon >0$, there exists an $L>0$ such that $x_j\in B(x_{\sssize L}, \epsilon), \forall j>L.$ With
respect to each $x_j$, there associates a point $(x_j,v_j)\in N.$ Thus $\{(x_j,v_j):j>L\}$ is a
sequence in the compact set $N\cap\overline {T^l B(x_{\sssize L},\epsilon)}$ and thus has a limit point 
$(x,v)\in N\cap\overline {T^l B(x_{\sssize L},\epsilon)}.$ Thus $x\in \pi (N)$ and $\pi (N)$ is closed.
\qed

\enddemo
The Lie group $I $ is a  subgroup of $G $ and $I_0=G_0$. We consider the group
$G /I_0 $ and examine the index of this coset space. 

\proclaim{Proposition 4.3}
The index of $I_0 $ in $G $ is  finite.
\endproclaim
\demo{Proof}
Let $\{g_j I_0 \}, g_j\in G $, be a sequence in the coset space $G /I_0. $ 
By Lemma 4.2,
$g_j(M)$  has non-empty intersection with any fiber $T^r_pM,  p\in M$.
Fix $p\in K$ and  take
a point $$q_j\in g_j(M)\cap  T^r_pM.\tag 4.2$$ Since $g_j^{-1}(q_j)\in M$ there exists a
$h_j^{-1}\in I_0$ such that  
$$h_j^{-1}\cdot g_j^{-1}(q_j):=a_j\in  K, \forall j.\tag 4.3$$ Let 
$f_j=g_j\cdot h_j\in G $ then $$f_j(a_j)=q_j\in T^r_pM.\tag 4.4$$

By the complete hyperbolicity proved in $\S 3$, $T^rM$ is a taut manifold which says that we can
extract a subsequence(we still call it $\{f_j\}$) that either converges uniformly on compact subsets of
$T^rM$ or diverges compactly, i.e., for any compact subsets $D_1,D_2$ of $T^rM$ there exists a large
$\Cal N$ such that
$f_j(D_1)\cap D_2=\emptyset, \forall j>\Cal N.$

Let $$diam (K):=\sup \{d_{\sssize K}(p,q): p,q\in K\}=\Cal R<\infty$$ denote the diameter of the
compact set $K$ with respect to the Kobayashi metric $d_{\sssize K}$ in $T^rM$. Then $K\subset
B_{\sssize K}(a_j,\Cal R),\forall a_j\in K.$ As each $f_j$ is an isometry of the Kobayashi
metric. We have $$ f_j( K)\subset f_j (B_{\sssize K}(a_j,\Cal
R))= B_{\sssize K}(f_j(a_j),\Cal R)=B_{\sssize K}(q_j,\Cal R), \forall
j.$$  Since $q_j\in T^r_pM$, it is therefore  possible to find a compact set $K'\subset M$
such that $$f_j( K)\subset T^r_{\sssize K'}M,\forall j,\tag 4.4    $$
where $T^r_{\sssize K'}M=\cup_{p\in K'} T^r_pM.$
Let $\Omega_k:=T^{r-\frac 1{k}}_{\sssize K'}(M).$ Suppose the sequence $\{f_j\}$ diverges compactly.
Then for any $k\in \Bbb N$ there exists an $\Cal N_k$ such that $f_j(K)\cap \bar \Omega_k=\emptyset$ for
all
$j>\Cal N_k.$ That is $f_j(K)\subset \{(x,v);x\in K', v\in T_xM, r-\frac 1{k}<|v|<r\}$ for all
$j>\Cal N_k$. Let
$k\to\infty$, then the sequence $\{f_j\}$ sends the compact set $K$ to smooth strictly pseudoconvex
boundary points.
By assumption,  the \gt\ $T^rM$
is not covered by the ball. The generalized Wong-Rosay
theorem in [K2] implies that  no subsequence of $\{f_j\}$ could diverge compactly. A contradiction.
Hence, by the tautness of $T^rM$, there
exists a subsequence of
$\{f_j\}$ converges uniformly to some
$f \in G $ in the topology of $G $. Hence, some subsequence of $\{g_j I_0 \}$ converges to
$g I_0 $ in the topology of $G /I_0 $. Thus, $G /I_0 $ is compact. Since $I$ and $G$ have the same
identity components, 
$I_0$ is open in $G$. The compactness  implies that  the index of $G /I_0 $ is finite.\qed
\enddemo
For any two anti-holomorphic involutions $\sigma$ and $\tau$ in $T^rM$. We  use the notation 
$(M,\sigma)$ to 
 indicate that the anti-holomorphic involution of the Grauert tube $T^rM$ with respect to the
center
$M$ is $\sigma$.
Following the argument in Proposition 7.3 of  [K2], we prove
that if  $\tau$  is another  anti-holomorphic involution of $T^rM$, then the least
positive integer $k$ such that  $(\sigma\cdot\tau)^k\in I $ is odd and in fact 
 $(\sigma\cdot\tau)^k=id$. 

Proposition 7.4 and Lemma 7.5 of [K2] work in the quasi-homogeneous case as well.  
There are at most a finite number of anti-holomorphic involutions
$\sigma_j$ in $T^rM$. For any  center
$(N,\sigma_j)$ of
$\Omega$, there exists $f\in G $ such that $f(M)=N$.
The $Isom (M)$ is isomorphic to  $Aut (T^rM) $ if and only if there is a unique
center $(M,\sigma)$ for $T^rM$. 

Given $x\in M$, we consider  the homogeneous submanifold $I_0\cdot x$ of $M$.
\proclaim{Lemma 4.4} For $x\in K$,  the orbit
$I_0\cdot x$ is a totally geodesic submanifold of $M$.
\endproclaim
\demo{Proof} We chose a slightly larger set $\hat K\subset M$ such that $K$ is relatively compact in 
$\hat K$.
Let $B(x,\epsilon)$ denote the intersection of $\hat K$ with the geodesic ball in $M$ centered at $x$
of radius
$\epsilon$. We may choose $\epsilon$ so small such that $B(x,\epsilon)$ is relatively compact in
$\hat K$ and there is no $g\in I_0$ sending a
point $y\in B(x,\epsilon)$ to a different point $w\in B(x,\epsilon)$.

By the fact that $I_0\cdot K=M$, it is clear that $I_0\cdot B(x,\epsilon):= U$ is an open subset of
$M$, hence is a totally geodesic submanifold of $M$ with the induced metric from $M$. Since there is
no
$g\in I_0$ sending a point $y\in B(x,\epsilon)$ to a different point $w\in B(x,\epsilon)$, we may
view $U$ as a product manifold of $I_0$ and $B(x,\epsilon)$. Using the respective induced metrics  in
$U$ and
$B(x,\epsilon)$ from the Riemannian metric of $M$, we are able to put a Riemannian metric on $I_0$
such that $U$ is the product manifold of $I_0$ and $B(x,\epsilon)$ with the product metric. 

Since $I_0\cdot x$ is a totally geodesic submanifold of the product manifold $I_0\cdot B(x,\epsilon)$
which is totally geodesic in $M$. Therefore, $I_0\cdot x$ is a totally geodesic submanifold of $M$.
\qed
\enddemo

Let's denote $M_x=I_0\cdot x.$
 As $M_x$ is a totally geodesic submanifold of $M$, the \gt\ $T^r(M_x)$ is a complex
submanifold of $T^rM$. Since $I_0$ acts transitively on $M_x$, 
the tangent space $T_z(T^rM_x)$ could be decomposed as, for any $z\in T^rM_x$,  $$T_z(T^rM_x)=T_z(I_0
\cdot z)+T_z(T_{\pi(z)}^r M_x).\tag 4.5$$ 
Following the   arguments in $\S7$ of [K-M]. We are able to construct a $G$-invariant strictly
plurisubharmonic  function $$\psi(z)=\sum_{j=1}^k \rho(g_j(z))$$ in $T^rM$ where
$\{g_1,\dots, g_k\}\in G$ so that $G/G_0=\{g_jG_0: j=1,\dots, k\}$ as shown in Prop. 4.3. 

Let $F_x:= \psi|_{T^r M_x}$ denote the restriction of
$\psi$ to $T^r M_x$ and  $\eta_x:= \psi|_{T_{\pi(z)}^r M_x}$  the restriction of
$\psi$ to ${T_{\pi(z)}^r M_x}.$ 
Since $\psi$ is constant in $I_0 \cdot z$, by the decomposition (4.5), a critical point of  $\eta_x$ is
a  critical point of $F_x$.

As $T^rM_x$ is a submanifold of $T^rM$, the restriction $F_x$ of $\psi$ to $T^rM_x$ is still 
 strictly plurisubharmonic. The real Hessian of  $F_x$ has at least $k=dim (I_0\cdot x)$ positive
eigenvalues and can be null on a subspace of $T_z(T^rM_x)$ of dimension at most $k$. Since $F_x$ is
constant on $T_z(I_0\cdot z),$ and $dim (I_0\cdot z)\ge k$, we have  $dim\ T_z(I_0\cdot z)=k$ and the
 real Hessian of $\eta_x$ is positive definite on the tangent space 
$T_z(T_{\pi(z)}^r M_x).$ By the Morse theory (c.f. p85, [C-E]) $\eta_x$ has at most one critical point.
Since $\eta_x$ is proper on the  fiber $T_{\pi(z)}^r M_x$,
 it follows that there is
exactly one critical point of $\eta_x$ which turns out to the the minimal point. Since $\eta_x\cdot
\sigma=\eta_x$, the minimum of $\eta_x$ occurs at $\pi (z).$ 

 On the other hand, the restriction $\psi_p$ of $\psi$ to the fiber $T^r_pM$ is proper, therefore
$\psi_p$ has local minimums in $T^r_pM$ which must also be critical points of $\eta_p$. Therefore,
each $\psi_p$ has exactly one minimal point at $p$. Let
$$L:= \max_{p\in K} \psi(p).$$
As $\psi$ is $G$-invariant, $M=I_0\cdot K$ and $N=f(M)$, we have $$\max_{p\in M}\psi(p)= \max_{p\in K}
\psi(p)=\max_{z\in N}\psi(z).$$

Let $\Cal A:=\{x\in K: \psi (x)=L\}.$ Pick $q\in \Cal A$  and $z\in N\cap
T_q^rM$. Since
$z\in N$,
$\psi (z)\le L$. On the other hand, $z\in T_q^rM$, $\psi (z)\ge\psi (q)=L.$ This can't be unless
$z=q$. That is, every maximal point lies in the intersection of different centers; 
 $q \in N\cap M$ and $q$ is a maximal point for $\psi|_{\sssize M}$.

By the quasi-homogeneity of $M$, the maximal points for $\psi|_{\sssize M}$ are $I_0\cdot \Cal A$ and
the maximal points for $\psi|_{\sssize N}$ are $f(I_0\cdot \Cal A)$. Therefore, $I_0\cdot \Cal A\subset
f(I_0\cdot \Cal A).$ Conversely, if we start from the center $N$, we will obtain $f(I_0\cdot \Cal
A)\subset I_0\cdot \Cal A$. Therefore $\psi|_{\sssize N}$ and $\psi|_{\sssize M}$ have the same
maximal point set.

 For any two centers $M$ and $N_j$,
there exists an
$g_j\in Aut (T^rM)$ such that
$g_j(M)=N_j.$ The above argument works for any center as well. Let $q\in \Cal A$, then $g_j(q)$ is a
maximal point for $\psi|_{\sssize N_j}$. Hence, $g_j(q)\in M, \forall j$ which implies 
$$\psi (q)=\sum_{j=1}^k \rho(g_j(q))=0.$$ 
 As $\psi$ is a non-negative function
and $q$ is a maximal point for $\psi|_{\sssize M}$. We conclude that $\psi|_{\sssize M}\equiv 0$.
Therefore, $g_j(M)=M$, for all $g_j\in Aut (T^rM). $

The following main theorem is thus proved.

 \proclaim{Theorem 4.5}
Let $M$ be  a \qh\ manifold. Let
$T^rM$ be the
\gt\ of radius
$r<r_{max}(M)$ which is not
 covered by the ball.  Then
$Isom(M)\approxeq Aut(T^rM)$ and $T^rM$ has a unique center.
\endproclaim

\enddocument